\documentclass[11pt]{article}
\usepackage{amsfonts, amssymb, amsmath, amscd, longtable, tabularx, multicol}
\usepackage{graphicx}
\newtheorem{theorem-intro}{Theorem}[section]

\newtheorem{corollary-intro}[theorem-intro]{Corollary}

\def\be{\begin{enumerate}}
\def\ee{\end{enumerate}}
\def\bi{\begin{itemize}}
\def\ei{\end{itemize}}
\def\xi{\mathcal X}

\title{Lorenz knots}
\author{Joan S. Birman, Department of Mathematics, Barnard-Columbia, 2990 Broadway \\
New York, NY 10027, USA, e-mail $<jb@math.columbia.edu>$}

\begin{document}
\maketitle

\begin{abstract}  
This is a review article on Lorenz knots.  It will appear in the Proceedings of a conference organized by Prof. Christophe Letellier, in honor of the $70^{\rm th}$ birthday of Professor Robert Gilmore, June 28-30, 2011, at the University of Rouen in  Rouen, France.  
\end{abstract}

\tableofcontents

%%%%%%%%%%%%%
%%%%%%%%%%%%%
%%%%%%%%%%%%%%
\section{Introduction}\label{JBS:introduction}  We begin, in $\S$\ref{JBSS:history}, with a very informal and intuitive review of the history of knot theory.  In $\S$\ref{JBSS:the search for invariants} we give precise definitions of knot and link type, and discuss some of the invariants  that were discovered over the years.  In $\S$\ref{JBSS:Special classes of knots and links}  we show how the link classification problem changes when the links in question are limited to special classes, with torus links as an example.  In $\S$\ref{JBSS:questions} this will to lead us naturally  to Lorenz links, the main subject of this review.  They are, in a very precise way, generalizations of torus links.

\subsection{History}\label{JBSS:history}At the end of the 19th century the well-known Scottish physicist Peter Guthrie Tait (1831-1901) had the idea that the periodic table might be explained by knotting in the `impenetrable ether'.  This lead him to study knots, and while he did not succeed in his original goal, he became fascinated  by the intricacies of knotting.  This was the beginning of the systematic study of knots, and (in a sense) a beginning of the part of topology that we call {\it Knot Theory}.  

 If we widen the discussion to include not just knots but also links,  there was even earlier work due to Gauss (1777-1855), who was interested in the question of how the current in a closed knotted wire affected that in another closed knotted wire when the two closed closed curves were `linked'.  A different but related question also appears in Gauss' work  \cite{JBEpple98}, where he computed the linking number of the earth's orbit with that of certain asteroids.  

We've given two reasons why knots and links were of interest in physics, but what about mathematics?   Mathematicians first became interested in knots because they are very easily visualized, and their complementary spaces $\mathbb S^3\setminus K$, where $K$ is a knot, are a rich source of examples of  3-manifolds.   While knot complements are non-compact, `Dehn surgery' on knot spaces, defined by associating a rational number to each  component of a link \cite{JBRo}, is a way to construct all closed 3-manifolds too.  We note that the phenomena of knotting occurs in every dimension (especially in codimension 2) as part of  a study  of the way that one manifold sits inside another.  

Knots also appear in other parts of mathematics where there is no obvious way, at least initially, to visualize knotting.   For example, an algebra  student who has learned  the definitions of a ring and a subring and has worked on some good examples  could probably understand that  if $R'$ is a subring of a ring $R$, then the inclusion map $R' \subset R$ might be interesting from a combinatorial viewpoint.  Type II$_1$ factors in operator algebras may be thought of analogues of  rings.  When the Jones polynomial was discovered in 1985,  its mathematical origins were in Operator Algebras.  It was later understood that it was crucial to the discovery that some sort of knotting-in-disguise was involved in the way that one type-II$_1$ factor sat inside another type-II$_1$ factor.  See Birman's article \cite{JBBir92} on the work of Jones for a fuller discussion.   

\begin{figure}[htpb!]
\centerline{\includegraphics[scale=.750]{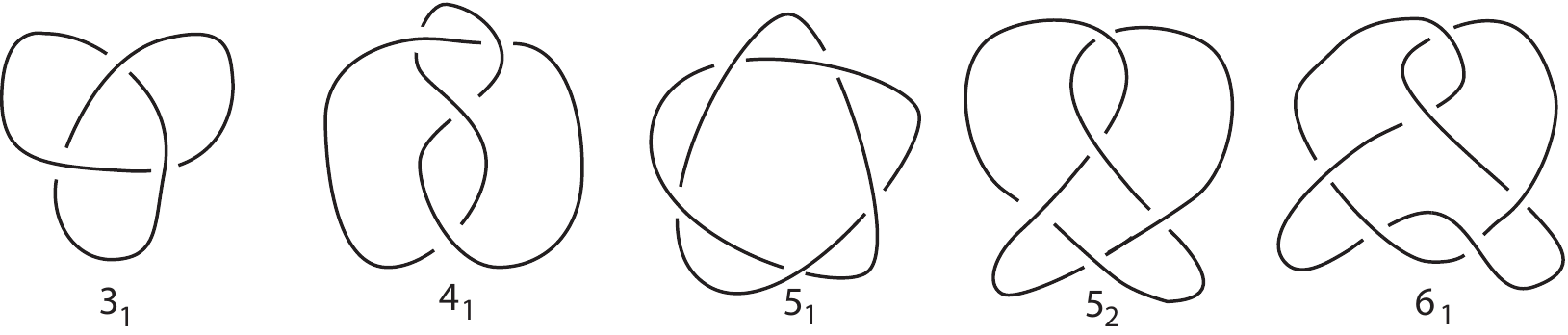}}
\caption{The first five knots in Tait's tables}
\label{JBknottables}
\end{figure}

Data from the early knot tables has survived to this day.  Figure~\ref{JBknottables} shows the first 5 examples from the 1893 knot tables.  Those tables have now been extended to `the first 1,701,936 knots', where the measure of complexity is the number of crossing points in a picture of the knot.  See \cite{JBTh85}.  For example, there is only one knot of crossing number 3, also only one with 4 crossings, but two with 5 crossings and so on.  
By `only one' is meant one (up to distinct pictures of the same knot).   This is a good moment to mention that, while minimum crossing number, over all possible projections, is by definition independent of the choice of a projection,   it has not turned out to be a very meaningful measure of complexity.  Indeed,  modern knot tables include ones where the measure of complexity reflects more subtle aspects of knot theory.  We will have more to say about that in $\S$\ref{JBss:hyperbolic volume} below.

\begin{figure}[htpb!]
\centerline{\includegraphics[scale=0.70] {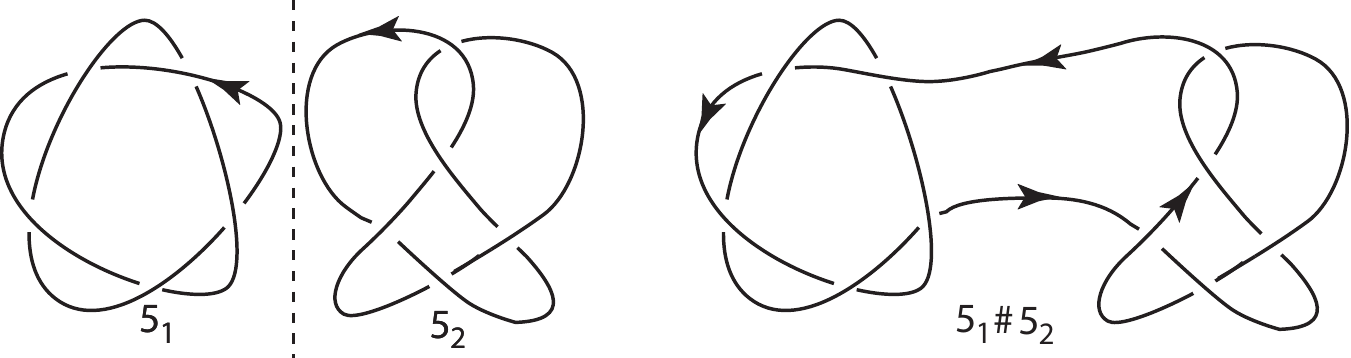}}
\caption{The knots $5_1$ and $5_2$ on the left are prime knots from Figure 1.  The right sketch shows the `composite' knot $5_1\#5_2$.  } 
 \label{JBcomposite}
\end{figure}

A knot is {\it composite} if you can obtain a representative by tying one
knot in a piece of string, and following it by another.    It's {\it prime}
if itÕs not composite.  See Figure~\ref{JBcomposite} for examples.
%%%
%%%
Since prime knots are the Òbuilding blocksÓ of all knots, that's why the knots in the tables  are all prime.  Note also that if $O$ is a picture of the `unknot', and $K$ is any fixed picture of a knot $K$, then $K\# O$ just gives us a new picture of $K$.  Caution: there do not exist knots $K, K'$, both different from the unknot, with $K \# K'$ the unknot $O$. That is,  knots form an addition semi-group, using the operation $\#$.    

Clearly the pictures in Figures~\ref{JBknottables} and \ref{JBcomposite} are not unique,  and equally clearly that problem will cause lots of trouble as the number of crossings increases, so it's time for a definition:
\begin{itemize}
\item    A {\it knot} $K$ is the image of a circle $S^1$ under an embedding  $e: S^1 \to \mathbb S^3$ or $\mathbb R^3$, that is $e(S^1) = K\subset \mathbb S^3$.   Two knots $K,K^\star$ have the same {\it type}  if there is a diffeomorphism of pairs $(K,\mathbb S^3) \to (K^\star,\mathbb S^3)$.   \end{itemize}
The term `knot' always means `knot type'.  However, we will be sloppy and speak of $K=e(S^1)$ when we really mean $K$ or any knot $K^\star = e^\star(S^1)$ which has the same knot type as $K$.   Our pictures represent some choice of a projection of the knot onto a plane,  with undercrossings and overcrossings distinguished.  

The problem of recognizing when two knot diagrams represent the same knot type is a highly non-trivial problem.    That's why we have a need for computable invariants that are independent of the projection.     
%%%%%%%%%%%%%%%
%%%%%%%%%%%%%%%
\subsection {The search for invariants}\label{JBSS:the search for invariants}
\begin{figure}[htpb!]
\centerline{\includegraphics[scale=.70]{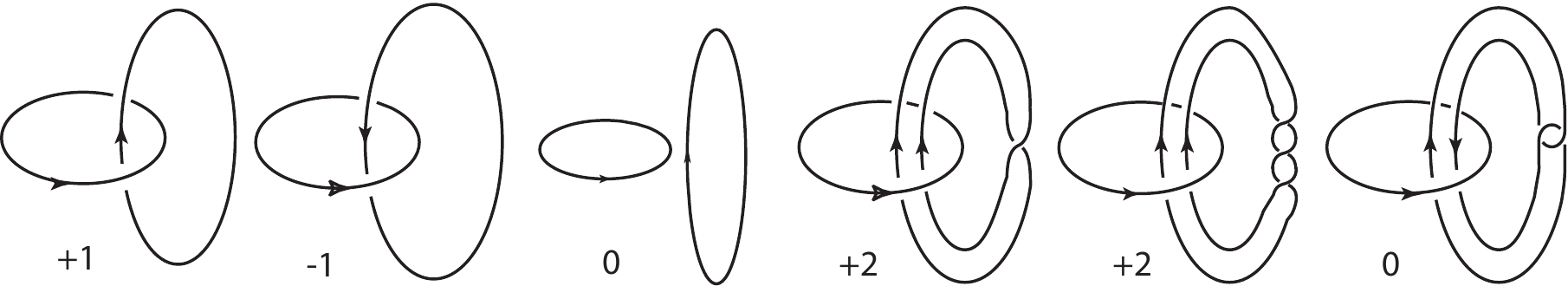}}
\caption{six 2-component links, with their linking numbers $\mathcal L(K,K')$.} \label{JBlinking-no}
\end{figure}

To define the earliest known invariant, that of Gauss, we need to generalize the concept of a knot:  
\begin{itemize}
\item [$\bullet$] A {\it link} $L$ is the image under an embedding of $\mu\geq 1$ disjoint circles in $\mathbb S^3$ or $\mathbb R^3$.   If $\mu=1$ it's a {\it knot}. 
Two links $L,L^\star$ have the same {\it type}  if there is a diffeomorphism of pairs $(L,\mathbb S^3) \to (L^\star,\mathbb S^3)$ or $(L,\mathbb R^3) \to (L^\star,\mathbb R^3)$.  
An {\it invariant} is any computable quantity which takes the same value on all representatives of the link type.
\end{itemize}
An example of a link-type-invariant is the number of components in a link $L$,  but Figure~\ref{JBlinking-no} shows six different 2-component links. 
The first subtle invariant is due to Gauss, who discovered an iterated line integral over 3-space $\mathbb R^3$ that is invariant under isotopy of the pair $(L,\mathbb R^3)$, where $L = K\cup K'$:  
% $$ \frac{1}{4} \int_K\int_{K'}  \frac{r_1 - r_2}{|r_1-r_2|^3} \cdot (dr_1 \times dr_2).$$
$$ \int_K\int_{K'} \frac{(x'-x)(dydz'-dzdy') + (y'-y)(dzdx'-dxdz') +
(z'-z)(dxdy'-dydx')}{((x'-x)^2+(y'-y)^2+(z'-z)^2)^{3/2}}$$
The Gauss integral turns out to be an integral multiple of $4\pi$, and the integer is known as the {\it linking number} $\mathcal L(K,K')$. The integers attached to the links in Figure~\ref{JBlinking-no} give $\mathcal L(K, K')$, but (as one might expect it's a useful but crude measure of the way that the two components link one-another.  Note that the linking number can be zero when the components are not visibly disjoint.
\begin{itemize}
\item [] {\bf Remark:} There is, in fact, an easy way to compute the Gauss linking number.  See Figure~\ref{JBsigns} for a way to assign $\pm 1$ to each crossings in a projection.  
Armed with this definition, color one component in each sketch in Figure~\ref{JBlinking-no} red, the other blue.  Using the sign convention that's given in Figure~\ref{JBsigns}, count  the algebraic sum of the signs at the crossings of red over blue (or blue over red, it will not matter), ignoring all monochromatic crossings.  It's the Gauss' linking number.   
 %%%%%%%%%%%%
 \begin{figure}[htpb!]
\centerline{\includegraphics[scale=.9]{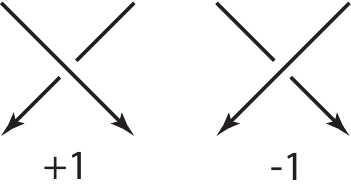}}
\caption{Sign conventions. }
\label{JBsigns}
\end{figure}
%%%%%%%%%%%%

{\bf Caution:} Our conventions come from Dynamical Systems, and are opposite to those used in Topology!  In particular, the value of the Gauss integral is based on conventions in topology.  Starting in $\S$2 of this paper, when we begin to discuss Lorenz links, we will adhere, consistently,  to the conventions in Figure~\ref{JBsigns}.
 \end{itemize}

Turning to knots, the {\it minimum crossing number} $c_{min}(K)$  over all possible projections $K$ of a knot is, by definition, a knot or link type invariant, however (except for very low crossing number) it turns out to be astonishingly difficult to compute, and to this day we do not have any real understanding of it in the general case, although if we place some restrictions (see $\S$\ref{JBSS:Special classes of knots and links} below) the situation changes.  

\begin{figure}[htpb!]
\centerline{\includegraphics[scale=.45]{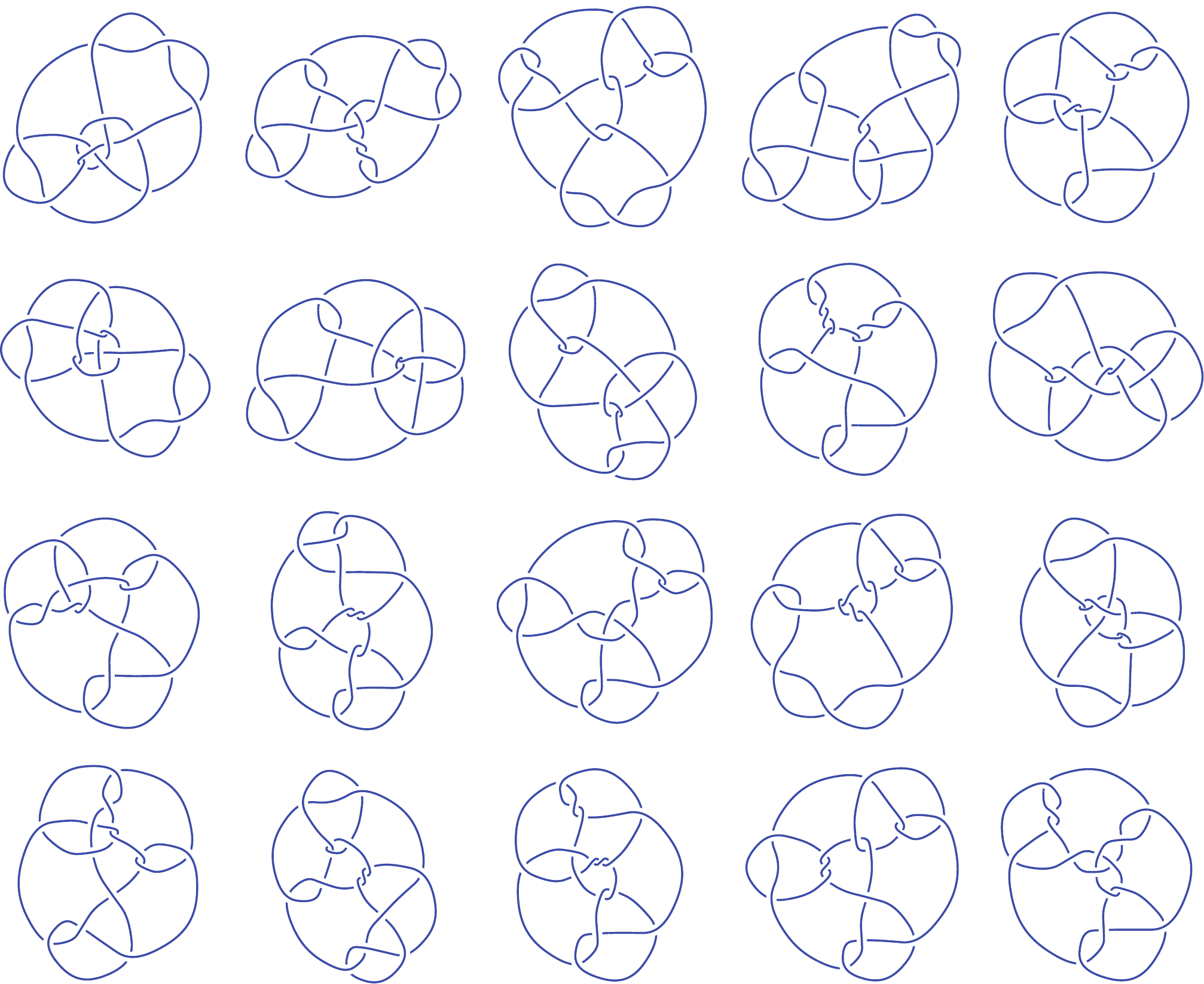}}
\caption{20 distinct knots with the same Jones polynomial (from M.Thistlethwaite)} 
\label{JBsameJonespoly}
\end{figure}

The Jones polynomial \cite{JBJo87} is a very sophisticated knot type invariant.  It's a Laurent polynomial with integer coefficients. Unlike the minimum crossing number,  it can be calculated from any projection, even though the calculation may be long and complicated if there are many crossings.   
Sadly, it's not a complete invariant.  See Figure~\ref{JBsameJonespoly}, which shows by example that distinct knots can have the same Jones polynomial.   
%%%
%%%

\subsection{Special classes of knots and links}\label{JBSS:Special classes of knots and links}  If we restrict our attention to special classes of knots or links, then the search for invariants can sometimes be very much simpler than in the general case. 

\begin{figure}[htpb!]
\centerline{\includegraphics[scale=.6]{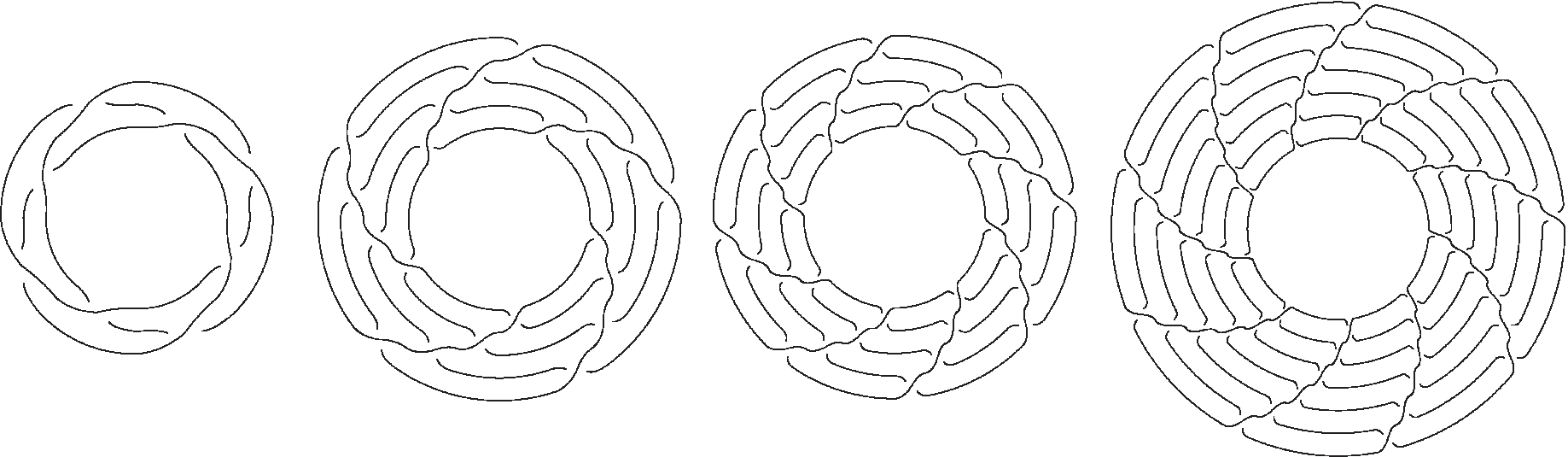}}
\caption{Four examples of torus knots and links. Their classifying  integer pairs, reading left to right, are $(3,5),(5,5),5,8),(7,9)$}
\label{JBtorusknots0}
\end{figure}

Torus knots and links are a special class.  They are defined to be knots and links that can be embedded on a standard torus of revolution in $\mathbb R^3$, as illustrated in Figure~\ref{JBtorusknots0}.  They are a well-understood class.   They are classified by a pair of integers $(p, q)$ (up to the indeterminacy $(p,q) \approx (q,p) \approx (-p, -q) \approx (-q, -p)$).    
%%%
%%%
The Jones polynomial, which we have already mentioned above,  of a $(p,q)$ torus link is interesting.  In the case when $p$ and $q$ are coprime, so that the link is a knot, it's given by the formula:
\begin{equation}\label{JBJones-poly-torus-links}
 V(t) = 
\frac{t^{ \frac{1}{2}(p-1)(q-1)}}
{1-t^2} 
(1 - t^{p-1} - t^{q-1} -t^{p+q})
\end{equation}  
The formula is considerably more complicated for links with more than one component, but for both knots and links it has a canonical form that's determined entirely by $p$ and $q$, so as far as the information that it contains it's a fancy way of saying $p$ and $q$.    But the point is that, unlike the integer pair $(p,q)$,  the Jones polynomial is defined for {\it all} knots and links.   

The {recognition problem} for torus links is to decide whether a knot or link is a torus link. It's a hard problem, but if you could solve it, then two integers would suffice to determine its type.  In general, the recognition problem for torus links is unsolved.  

{\it Positive or negative knots and links} are knots and links that admit a projection in which every crossing is consistently signed.  Like torus links (which are either positive or negative, depending on the defining parameters $p$ and $q$),  they are a proper subset of all knots and links, and like torus knots and links their study has been very fruitful.   It would be a bit of a diversion to discuss their special properties right now, so we ask the reader to be patient (or to look ahead to $\S$ \ref{JBapplications}      below,) where we shall have more to say about them.

%%%%%%%%%%%%
%%%%%%%%%%%%
\subsection{Some natural questions} \label{JBSS:questions}
When we encounter a new and interesting family of knots and links, especially one which arises in several different, seemingly unrelated settings (as we shall see, Lorenz knots and links are such a class), some questions we would like to answer are:  

\begin{itemize}
\item Which knots and links occur?  
\item  Is there a solution to the recognition problem for the class?
\item Does the knotting have meaning as regards the setting in which they were discovered? 
\item Is the mathematics of this class of knots related to phenomena being studied in other contexts?
\end{itemize}
We will address these questions, as they relate to the class of knots and  links that are the main subject of this review, `Lorenz'  knots and links,  in the pages that follow.

%%%%%%%%%%%%%%%
%%%%%%%%%%%%%%%
%%%%%%%%%%%%%%%
\section {Introducing Lorenz knots and links.} \label{JBSchaos}
In this section we introduce  Lorenz knots and links and prove that they include all torus knots and links.

%%%%%%%%%%%%%%
%%%%%%%%%%%%%%
\subsection{A chaotic flow on 3-space}
In a famous paper written in 1963 the meteorologist E. N. Lorenz asked the question: ``Is weather fundamentally  deterministic?"  It was clear that there were many many variables that affected the weather, but if we knew all of them, and the precise equations that governed weather patterns, would it really be {\it predictable} in a precise way?   Giving the question a slightly different twist, does {\it deterministic} $\Longrightarrow$ {\it ultimately periodic}? 
If so, then weather can't be deterministic.   He set out to investigate this question.  See his very readable article \cite{JBLor63}. 
 
His starting point was the Navier-Stokes equations, which he knew described the dynamics of a viscous, incompressible fluid.  
It had been used to model weather, ocean currents, water flow in a pipe, flow around an airfoil, motion of stars inside a galaxy.   All of those are very complex problems.  His idea was this: if a problem is too complicated to study it profitably, try to simplify it, preserving its essential features. discarding those which are unimportant.    Lorenz was a `theoretical meteorologist', and in this instance he was thinking like a mathematician, feeling free to modify the problem until it turned into a problem that might be amenable to study.
So he discarded variables, modified the equations, testing each time he changed things by using numerical integration.  After a very long (and very interesting) search he found a very simple system with the key features that he had wanted to retain.  
 His modified  equations turned out to govern fluid convection in a very thin disc that was heated from below and cooled from above.  The modified equations were a system of 3 ordinary differential equations in 3-space $\mathbb R^3$, and time $t$:
 \begin{equation} \label{JBLorenz-equations}
\frac{dx}{dt} = 10(y-x), \ \ \frac{dy}{dt} = 28x-y-xz, \ \ \frac{dz}{dt} = xy - \frac{8}{3}z  
\end{equation}
\noindent  There are three space variables: $x,y,z\in \mathbb R$ and their meaning is:

\qquad  $x$ = vertical temperature variation.

\qquad $y = $ horizontal temperature variation.

\qquad  $z$ = rate of convective overturning.

\noindent  The numbers 10, 20 and 8/3 were particular choices of parameters.    The equations are very robust, that is one may vary the parameters  10, 28,  and 8/3  in an open set without changing the features that he wanted to study in the solutions. 

\begin{figure}[htpb!]
\centerline{\includegraphics[angle=-70,scale=.25]{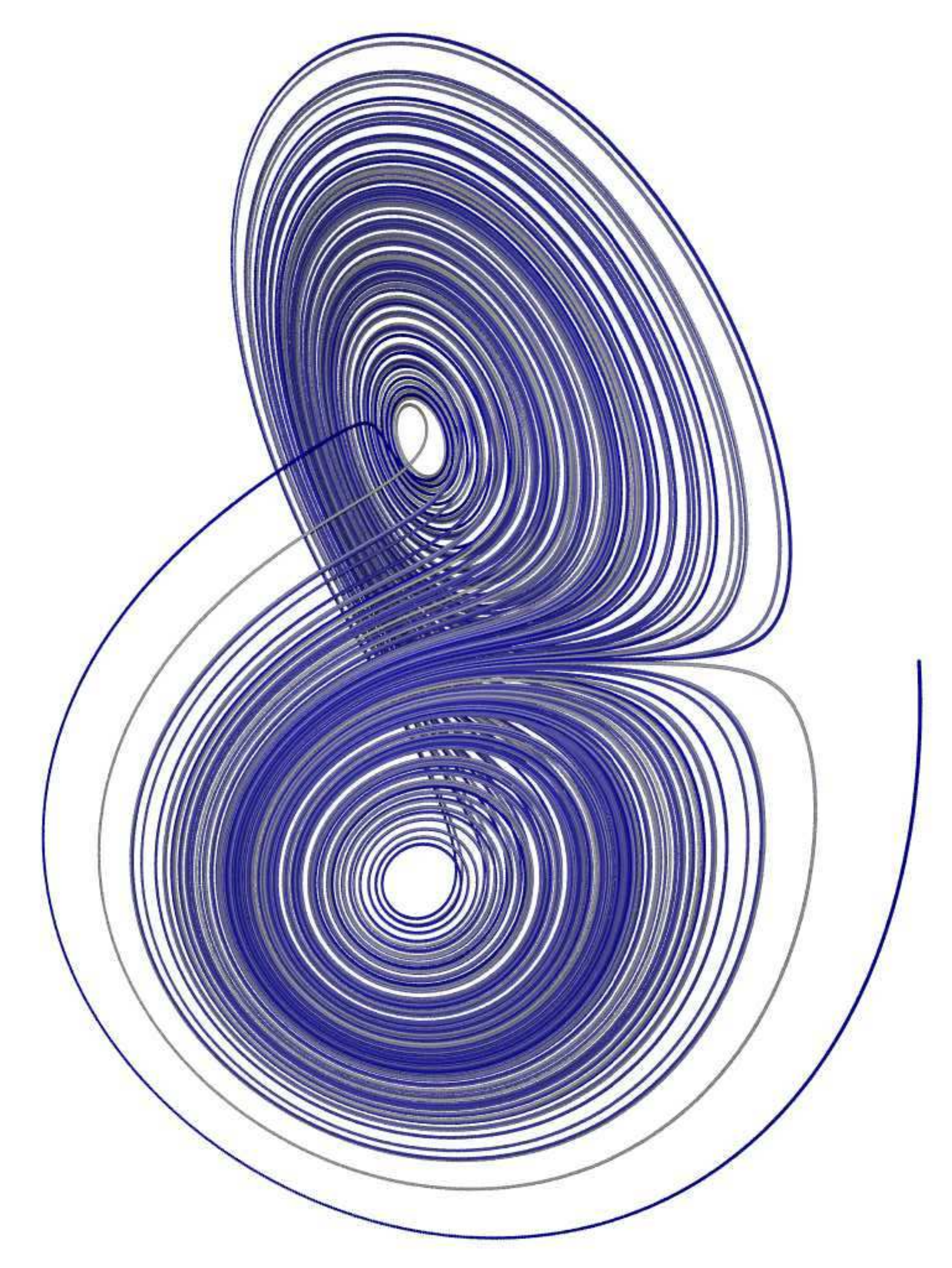}}
\caption{A beautiful picture of the orbits in the flow on $\mathbb{R}^3$ determined by
the Lorenz equations}\label{JBbeautifulpic}
\end{figure}

Of course, one can integrate a system of ODE's in $\mathbb R^3$, and if you do the orbits (with time as a parameter) determine a flow.  A {\it Lorenz knot} (resp. {\it link}) is defined to be a closed periodic orbit (resp. finite collection of orbits) in the flow $\lambda^t:\mathbb R^3\to \mathbb R^3$ determined by equations (\ref{JBLorenz-equations}).   Figure~\ref{JBbeautifulpic}  is a picture of what Lorenz saw when he integrated the equations (\ref{JBLorenz-equations}).  

Lorenz proved that there exists a bounded region $\mathcal A \subset \mathbb R^3$ such that forward trajectories ultimately enter $ \mathcal A$, and once they have entered it they stay there.  The region $\mathcal A$ is a neighborhood of the butterfly-shaped family of orbits that we see in Figure~\ref{JBbeautifulpic}.  Notice the layering of orbits.  
The subset $\mathcal A \subset \mathbb R^3$ is an {\it attractor.}    Now comes a key point:
orbits are extremely sensitive to initial conditions.  If you start at two points that are close to one-another,  their positions may be very far apart at a later time, even though both are inside $\mathcal A$.  Indeed, the Lorenz flow has become a prototype for a `chaotic flow'.  

\begin{figure}[htpb!]
\centerline{\includegraphics[scale=.8]{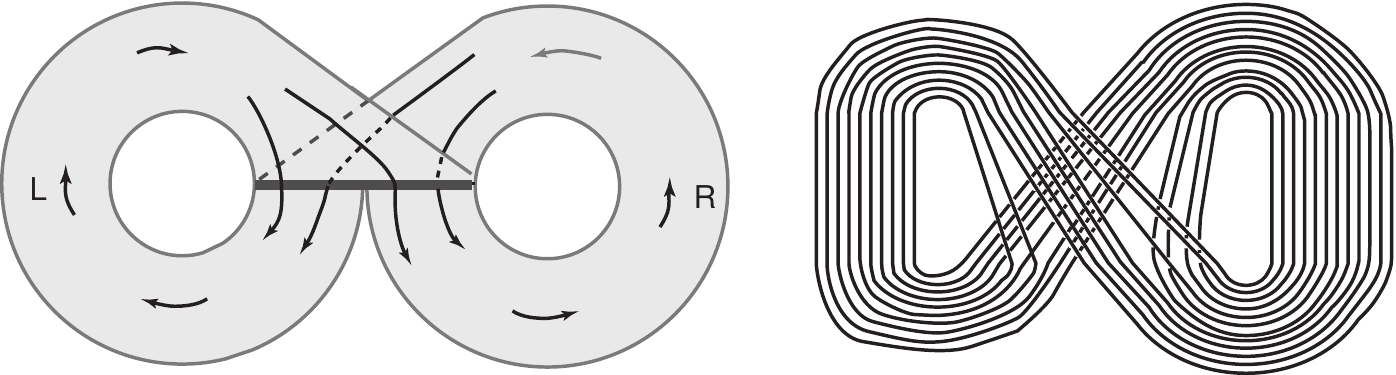}}
\caption{The left sketch shows the Lorenz template, and the right sketch shows a typical collection of orbits. The orbits are to be compared to those in Figure~\ref{JBbeautifulpic}}
 \label{JBorbits-suggest-template}
 \end{figure}

A  {\it template} is a branched 2-manifold with boundary, which is embedded in 3-space $\mathbb R^3$.  See \cite{JBGhy07} for an introduction and 
thorough review of the literature on this topic.  Figure~\ref{JBorbits-suggest-template} 
shows the {\it Lorenz template}.  One can almost see the template in Figure~\ref{JBbeautifulpic}.  The knots that are depicted to the right of the template in Figure~\ref{JBorbits-suggest-template} are typical orbits supported on the template in the left sketch.   All crossings are the same sign, and we will, from now on, adhere to the conventions in Figure~\ref{JBsigns}, i.e. all crossings are {\it positive}.  

 The concept of a template first appeared in the early work of Guckenheimer and Williams,  in the 1970's, in their studies of the Lorenz flow.    They constructed an embedded branched surface in $\mathbb{R}^3$ that  supports a `semi-flow', that is the flow on the branched 2-manifold is oriented and the orientation is not reversible. They reasoned, with very careful estimates and numerical data to back up their ideas, that every finite subset of the closed orbits in the Lorenz flow flow ought to project, simultaneously and disjointly, onto the template.  It took over 30 years before indirect mathematical proofs of the existence of the template appeared,  in the work of  Tucker\cite{JBTu02} and of Ghys \cite{JBGhy07}.    Using their results,  a Lorenz knot can be defined to be a simple closed curve in 3-space that embeds in the Lorenz template, and a Lorenz link is a finite collection of disjoint simple closed curves carried by the template.  From that viewpoint, we can think of them as knots and links in 3-space, deforming them as we wish, as long as the isotopy extends to an isotopy of 3-space.  We can also bring all the machinery of knot theory to bear on their study. 

%%%%%%%%%%%%%
%%%%%%%%%%%%%
\subsection{Every torus link is a Lorenz link}  \label{JBtorus links are Lorenz}
\begin{figure}[htpb!]
\centerline{\includegraphics[scale=.60]{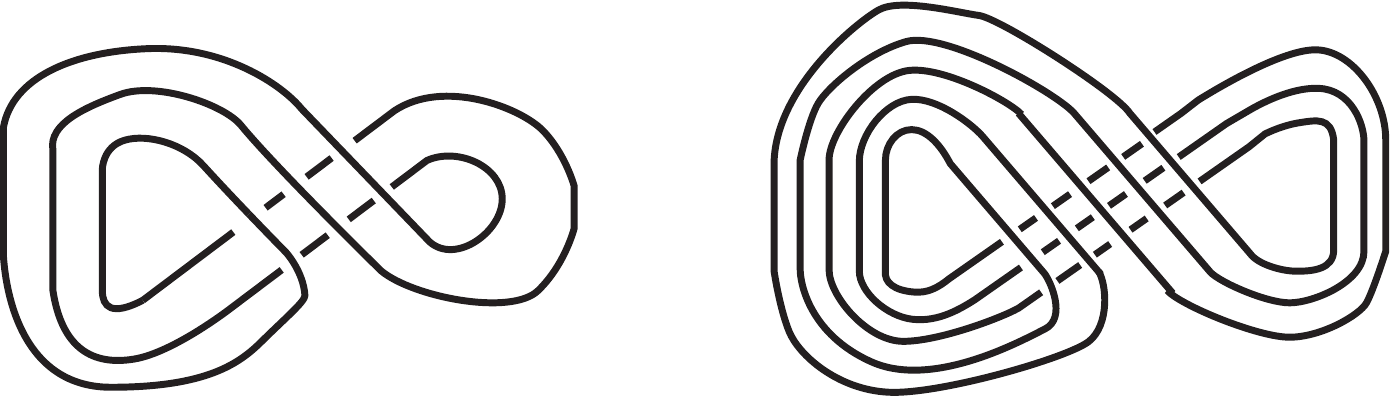}}
\caption{New pictures  torus knots of {\it type (2,3)} (on the left) and {\it type (3,5)} (resp. right)} 
\label{JBtorusknots1}
\end{figure}

We observed, earlier, that torus knots and links are a very special class.  It's a little bit of a surprise to learn that in fact every torus knot is a Lorenz knot!
The proof is surprisingly easy. It begins in Figure~\ref{JBtorusknots1} with a picture of two knots which are clearly Lorenz, because they can be embeded on the template. 
These are not the familiar projections of torus knots that we saw earlier, so how can we recognize them?  The sketches in Figure~\ref{JBtorusknots2}
shows us how to do it:
\begin{figure}[htpb!]
\centerline{\includegraphics[scale=.80]{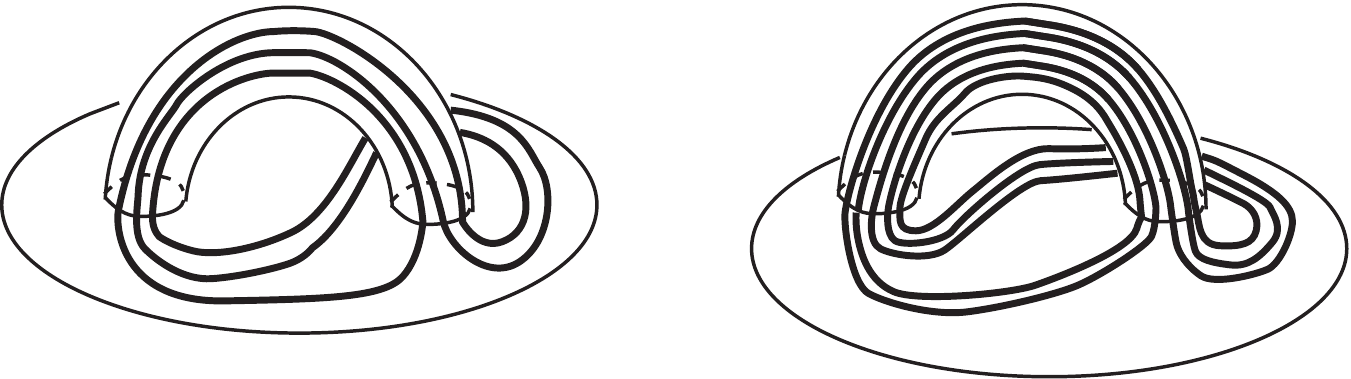}}
\caption{The  projection of Figure~\ref{JBtorusknots1} suggests, immediately, a way to embed these particular torus knots on a torus (minus a disc). Cap the single boundary component with a disc to get a closed torus.} 
\label{JBtorusknots2}
\end{figure}
Since the picture obviously generalizes from (2,3) and (3,5) to every integer pair $(p,q)$ we have learned that the set of all Lorenz knots and links is at least as big as the set of all torus knots and links.

%%%%%%%%%%%%%%
%%%%%%%%%%%%%%
%%%%%%%%%%%%%%
\section{Parametrizing Lorenz knots and links}   We describe three ways to parametrize Lorenz knots.  The order in which we choose to present them relates to our own taste, and does not reflect the order in which they were discovered historically.
\subsection{The first parametrization, Lorenz braids} \label{JBfirst parametrization}

\begin{figure}[htpb!]
\centerline{\includegraphics[scale=.60] {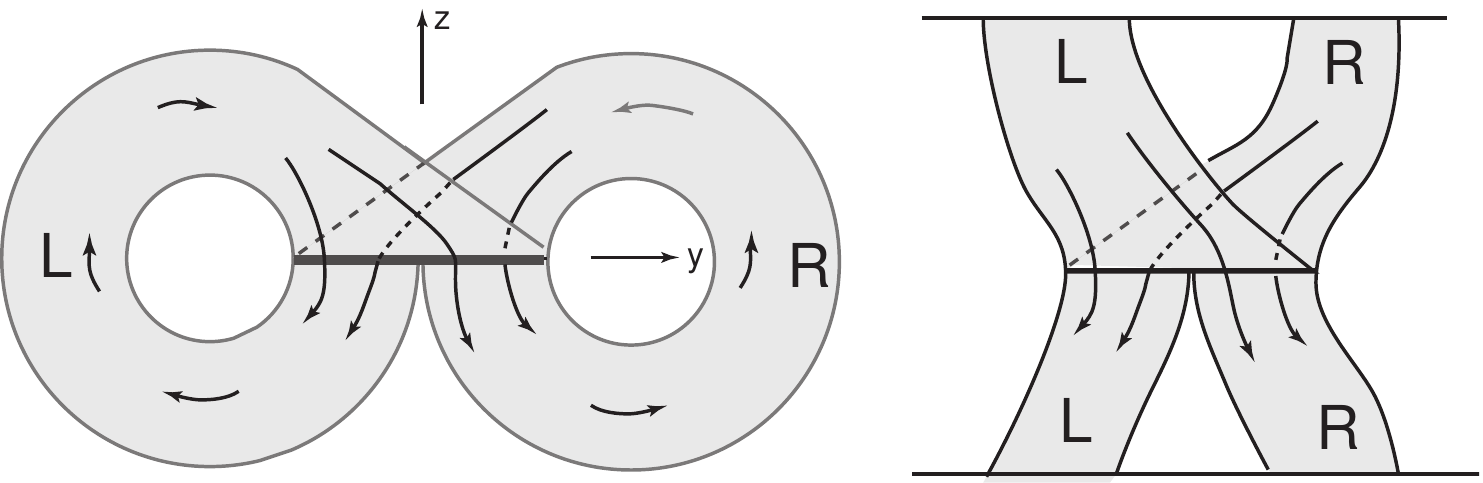}}
\caption{The left sketch is the Lorenz template.  The right sketch shows it cut-open to the Lorenz braid template.}
\label{JBbraidtemplate}
\end{figure}

The first paper \cite{JBBiWil83} in which Lorenz knots and links were studied as a class was written by Birman and Williams, and appeared in 1983, and the parametrization of Lorenz links by Lorenz braids was introduced there.  
In that paper a Lorenz link was defined to be any finite collection of closed orbits
on the Lorenz template, which supports a semiflow.   We now introduce a closely related branched 2-manifold, with boundary, which we call the {\it Lorenz braid template}.  
It can be seen in the right sketch in Figure~\ref{JBbraidtemplate}.  
%%%
%%%

Starting at the top of the braid template, the left and right branches of the template are stretched, and overlap along a horizontal {\it branch line},  then split apart, with one branch on the left and the other at the right.  The four `bands' are labeled L and R (at the top) and L and R again (at the bottom).   
A {\it Lorenz braid} is any finite set of braid strands that embeds on the Lorenz braid template.  

\begin{figure}[htpb!]
\centerline{\includegraphics[scale=.40] {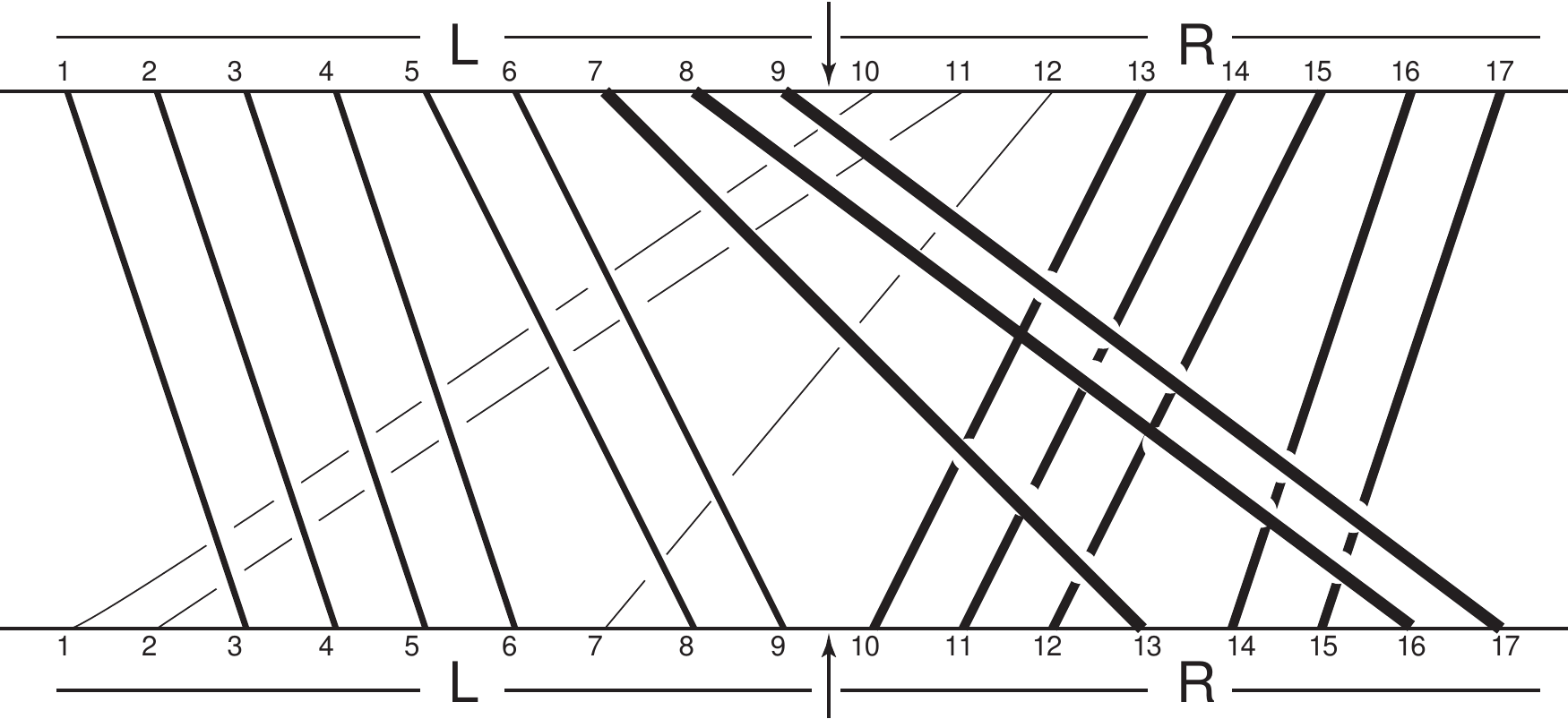}}
\caption{A typical Lorenz braid. }
\label{JBLorenzbraid}
\end{figure}

An example is given  in Figure~\ref{JBLorenzbraid}.  In the example there are 6,3,3,5  strands of type LL, LR, RL, RR.

%The important features of Lorenz braids are:  
%\begin{itemize}
%\item  Strands are oriented top to bottom, as dictated by the Lorenz flow.  
%\item All crossings are positive and two strands cross at most once, with all crossings being between a single overcrossing strand and a single undercrossing strand.  This implies that the braid is determined by the permutation associated to its braid strands.
%\item  A Lorenz braids can be parametrized by a sequence of positive integers $((p_1,q_1),(p_2,q_2).\dots,(p_r,q_r))$, where each $p_i<p_{i+1}$.  These parameters have the following meaning: the $i^{\rm th}$ group of overcrossing strands contains $q_i$ parallel strands, and if strand $j$ is in this group, then strand $j$ begins at position $j$ and ends at $j + p_i$.
%\item We make the following assumption: the first overcrossing strand ends at position $\geq 2$ and the final overcrossing strand ends at the last braid point. 
%\end{itemize}
%These are necessary and sufficient conditions, that is a braid with these properties closes to a Lorenz knot or  link, and every Lorenz link is the closure of such a braid. 

We assume that the projection, determined by the ODE's in (\ref{JBLorenz-equations}) and illustrated  in Figure~\ref{JBLorenzbraid}, is onto the region in the $xz$ plane between $z=0$( at the top) and $z=-1$) (at the bottom), and that the initial points of the braid strands have $x$ coordinates $1,2,3,\dots$.  Every crossing is therefore between a single overcrossing strand and a single undercrossing strand.  Those crossings are always positive, where we follow the sign convention that was given in Figure~\ref{JBsigns}.  Note that on each overcrossing
strand the final position will always be bigger than the  initial position, therefore if the $i^{th}$ overcrossing braid strand starts at (say) $x_i=i$, it will end at $x'_i > x_i$.  

Let $q_j$ be the number of strands in the
$j^{\rm{th}}$ group,  and let $p_j$ be the difference between the $x$-coordinates at the start of any one strand in the group 
at the end.   Suppose that there are $k$ such groups in the braid.  Then the numbers $(p_j,q_j), \ j=1,\dots,k$ completely determine the braid.  For this reason we can parametrize the braids by the $2k$ integers 
$ ( p_1,q_1), \dots, (p_k,q_k)$, where each $q_i\geq 1$, also  each $p_i < p_{i+1},$  and finally (a convenient assumptions that eliminates trivial cases)  $p_1\geq 2$ and $q_k\geq 2$.   
Here $(p_i,q_i)$ means that $p_i,\ldots, p_i$ repeated $q_i$ times.  This is our {\it first parametrization} of Lorenz links. In the example, the parameters are $((2,4),(3,2),(6,1),(8,2))$.

Observe that in the special case $k=1$ this reduces to a single integer pair $(p,q)$, and the associated closed braid determines a torus link, as we had already shown in $\S \ref{JBtorus links are Lorenz}$.  Thus Lorenz links are a generalization of torus links!  We shall have more to say about this in 
$\S \ref{JBthird parametrization}$.

\subsection{The second parametrization:symbolic dynamics and LR words} \label{JBsecond parametrization}
\begin{figure}[htpb!]
\centerline{\includegraphics[scale=.5] {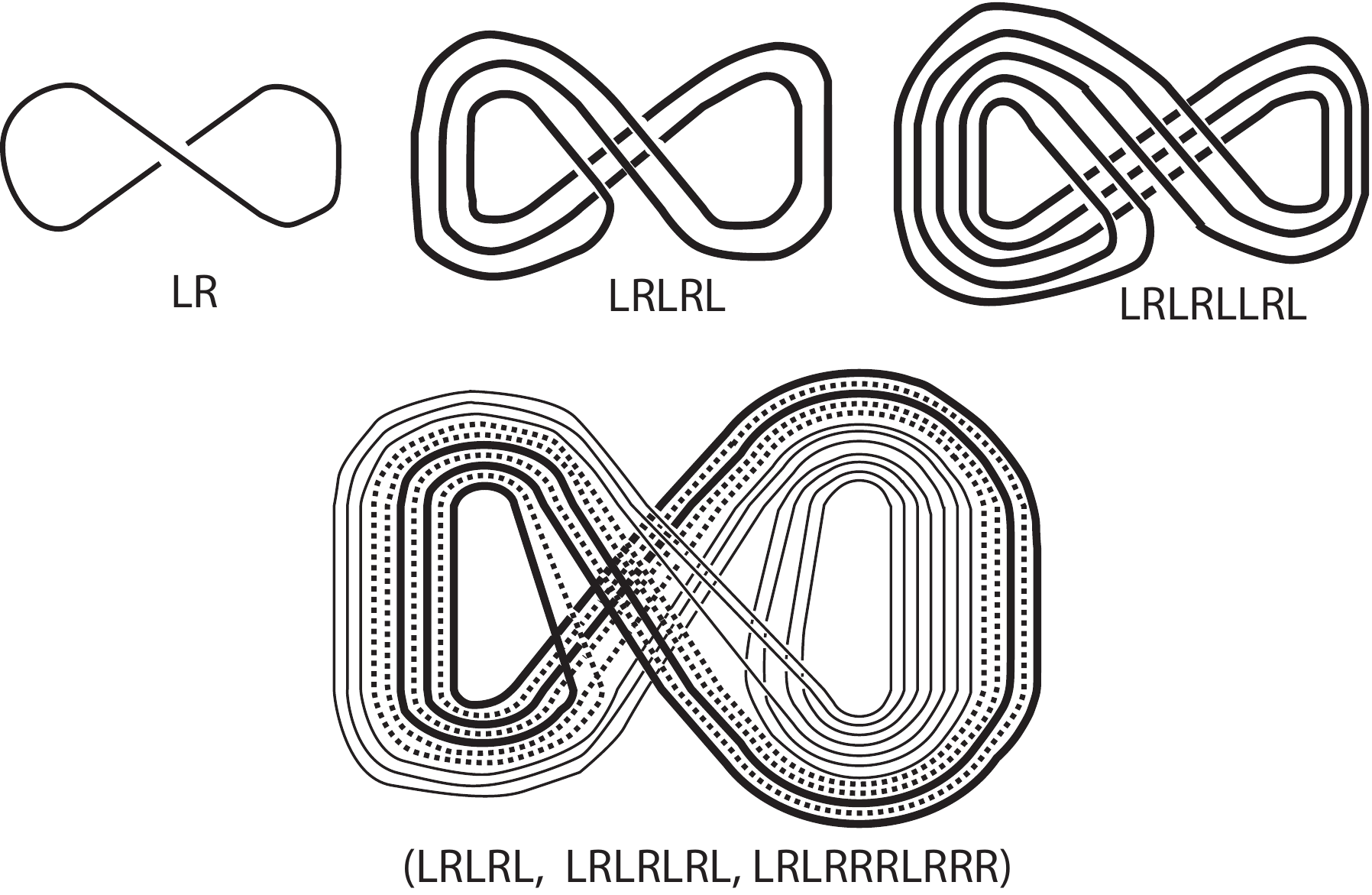}}
\caption{Examples of cyclic LR words and orbits in the Lorenz template.}
\label{JBwords1}
\end{figure}

The Lorenz template (see the left sketch in Figure~\ref{JBorbits-suggest-template})  tells us, immediately, that we may associate a word or words in the symbols $L$ (for left) and $R$ (for right) to any a closed orbit in the Lorenz flow.  There is no natural starting point for an orbit, therefore the $LR$ words used to describe it are cyclic.  See Figure~\ref{JBwords1}.  
This is an example of {\it symbolic dynamics},  a tool introduced to study solutions to ODE's like the Lorenz equations, which can't be integrated in closed form.   The bottom sketch in Figure~\ref{JBwords1} shows a link with 3 components, each described by its word.  The thickest line is  U = LRLRL  (that was our trefoil).  The dotted line: is V = LRLRLRL (a type (3,5) torus knot).   The thinnest line is   W = LRLRRRLRRR  (a type (-2,3,7) `pretzel knot').
It's clear that we don't want our word to be periodic.   As for links, we don't want the word for one component to be a power of that for another. 

Are there other restrictions?  The answer is `no'.   In the 1970's R. Williams proved \cite{JBBiWil83} that the correspondence between left-right words and closed orbits is a very strong one. He proved that a family  $\{W_1,\dots,W_k\}$ of cyclic words represents a Lorenz link if and only if no  $W_i$ is periodic, and no $W_i$ is (up to cyclic permutation) a power of any $W_j$.     

\begin{figure}[htpb!]
\centerline{\includegraphics[scale=.50] {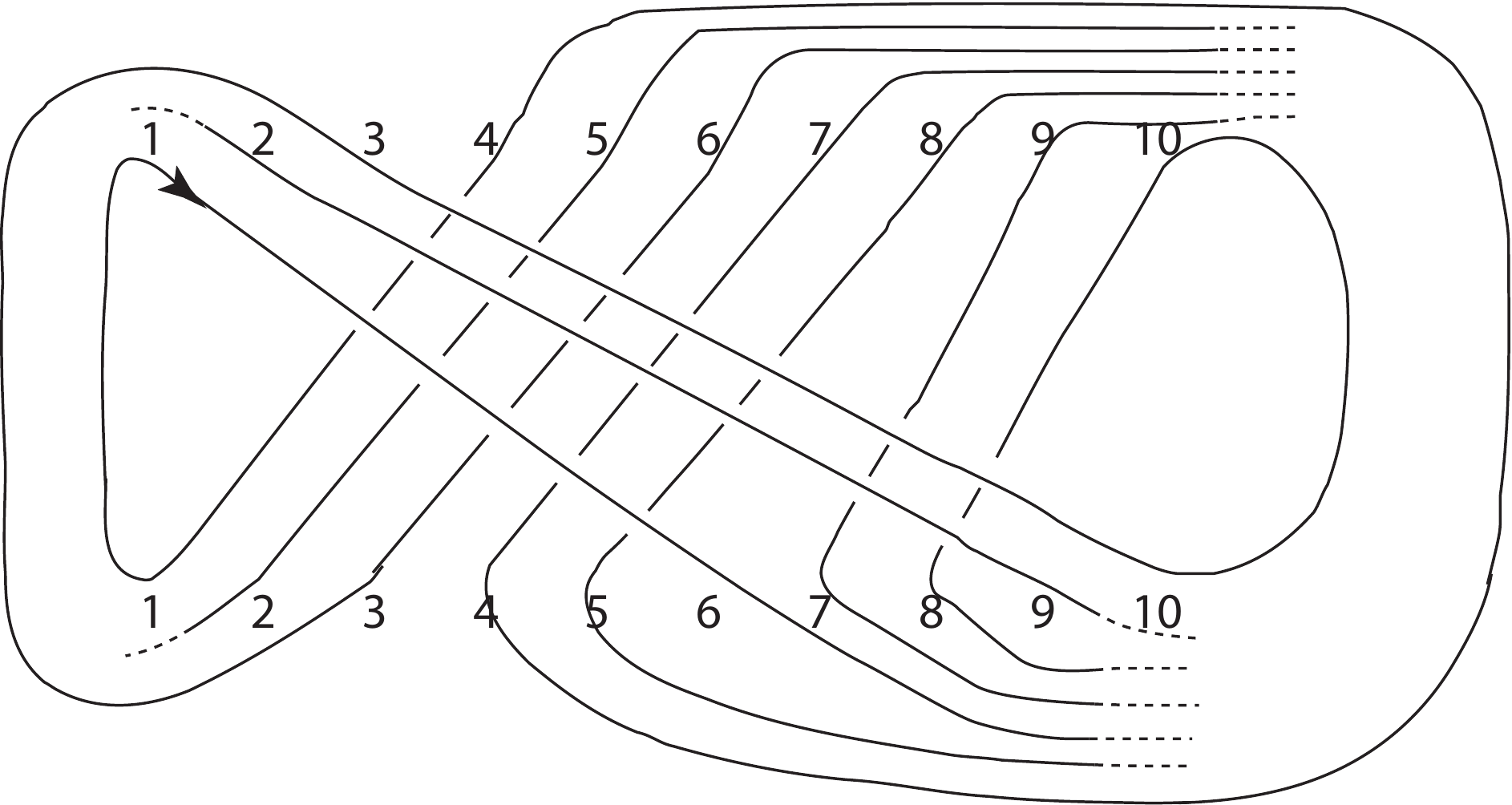}}
\caption{The Lorenz knot associated to $U$ = LRLRRRLRRR.}
\label{JBwords2}
\end{figure}

We suggest his proof, in Figure~\ref{JBwords2}, by showing how to recover the orbit from the word, in the example $U = $LRLRRRLRRR.   
To understand how we constructed the knot, the word $U$ = LRLRRRLRRR determines its 10 cyclic permutations $U = U_1,U_2,\dots,U_{10}$ , in their natural order.  Reorder these 10 words lexicographically, using the rule L$< $R.  The new position $\mu_i$  is given after each $U_i$: \\
 $U_1$ =  LRLRRRLRRR \ \ 1    \hspace {0.8 in}   $U_6$ =  RLRRRLRLRR \ \ 5 \\
$U_2$ =  RLRRRLRRRL \ \ 6      \hspace {0.8 in}  $U_7$ =  LRRRLRLRRR   \ \ 2\\
$U_3$ =  LRRRLRRRLR\ \ 3    \hspace {0.81 in}  $U_8$ =  RRRLRLRRRL \ \  9 \\
$U_4$ =  RRRLRRRLRL \ \ 10     \hspace {0.7 in}  $U_9$ =  RRLRLRRRLR \ \ 7 \\
$U_5$ =  RRLRRRLRLR  \ \ 8   \hspace {0.8 in} $U_{10}$ =  RLRLRRRLRR \ \ 4 \\
This determines a new cyclic order (1,6,3,10,8,5,2,9,7,4),  where strand $\mu_i$ begins at $\mu_i$ and ends at $\mu_{i+1}$.   We have constructed a {\it permutation braid}. 
A braid strand $i$ is an {\it overcrossing} strand if and only if  $\mu_i < \mu_{i+1}$, otherwise it's an undercrossing strand.  There are 3 (resp. 7) over (resp. under) -crossing strands in the example, hence 3 strands that go around the left ear and 7 around the right ear.  Starting with the permutation braid and traversing its associated closed braid, we recover the cyclic word LRLRRRLRRR.

\subsection{The third parametrization:  positive twisted torus links} \label{JBthird parametrization}

To explain the third parametrization, we return to our favorite example, the one in Figure~\ref {JBLorenzbraid}. There are 4 types of strands, those that begin at the left (or right) and those that end on the left (or right), and this will always be the case.  We label the four types of strands type LL, LR, RR and RL.

\begin{figure}[htpb!]
\centerline{\includegraphics[scale=.60] {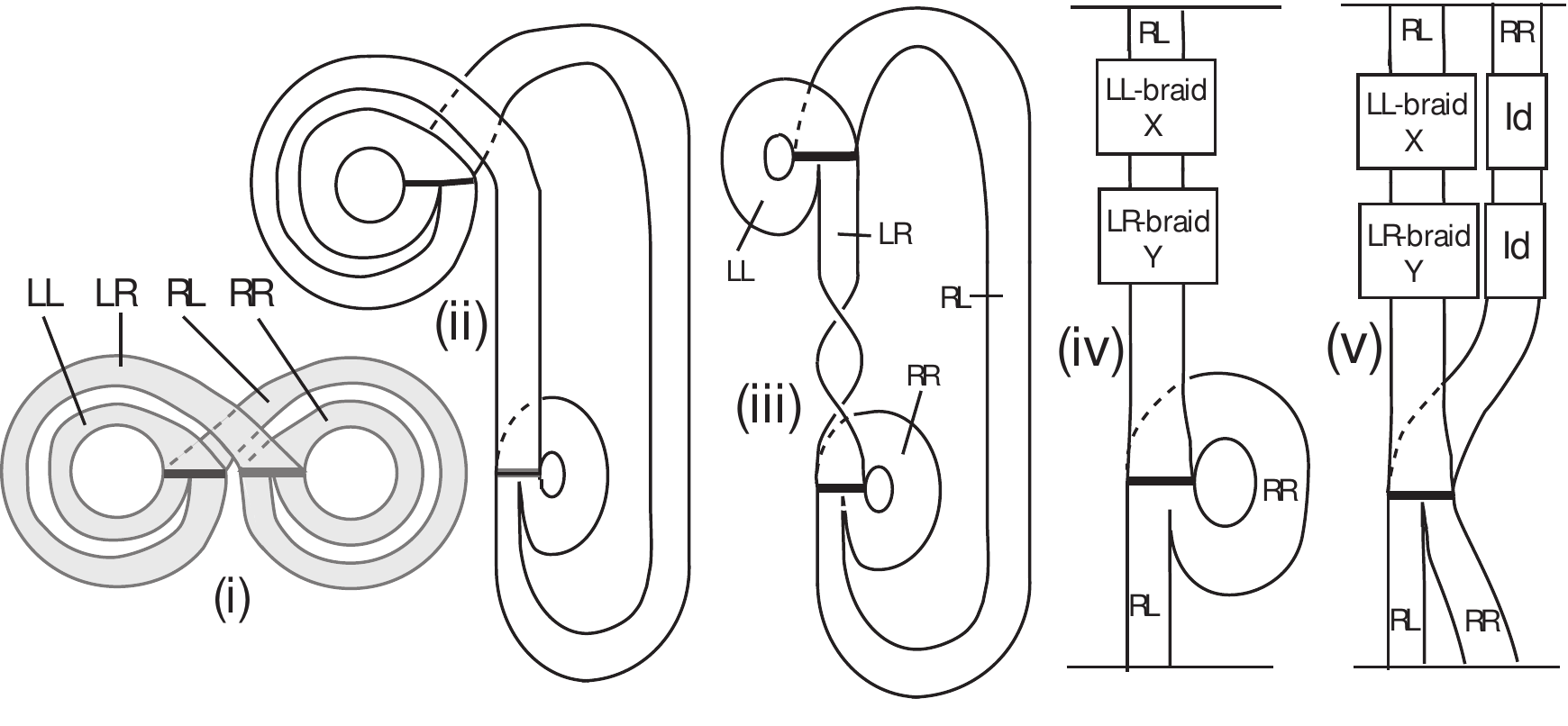}}
\caption{The steps in cutting open the Lorenz template} 
\label{JBcutting-open}
\end{figure}

\begin{figure}[htpb!]
\centerline{\includegraphics[scale=.60] {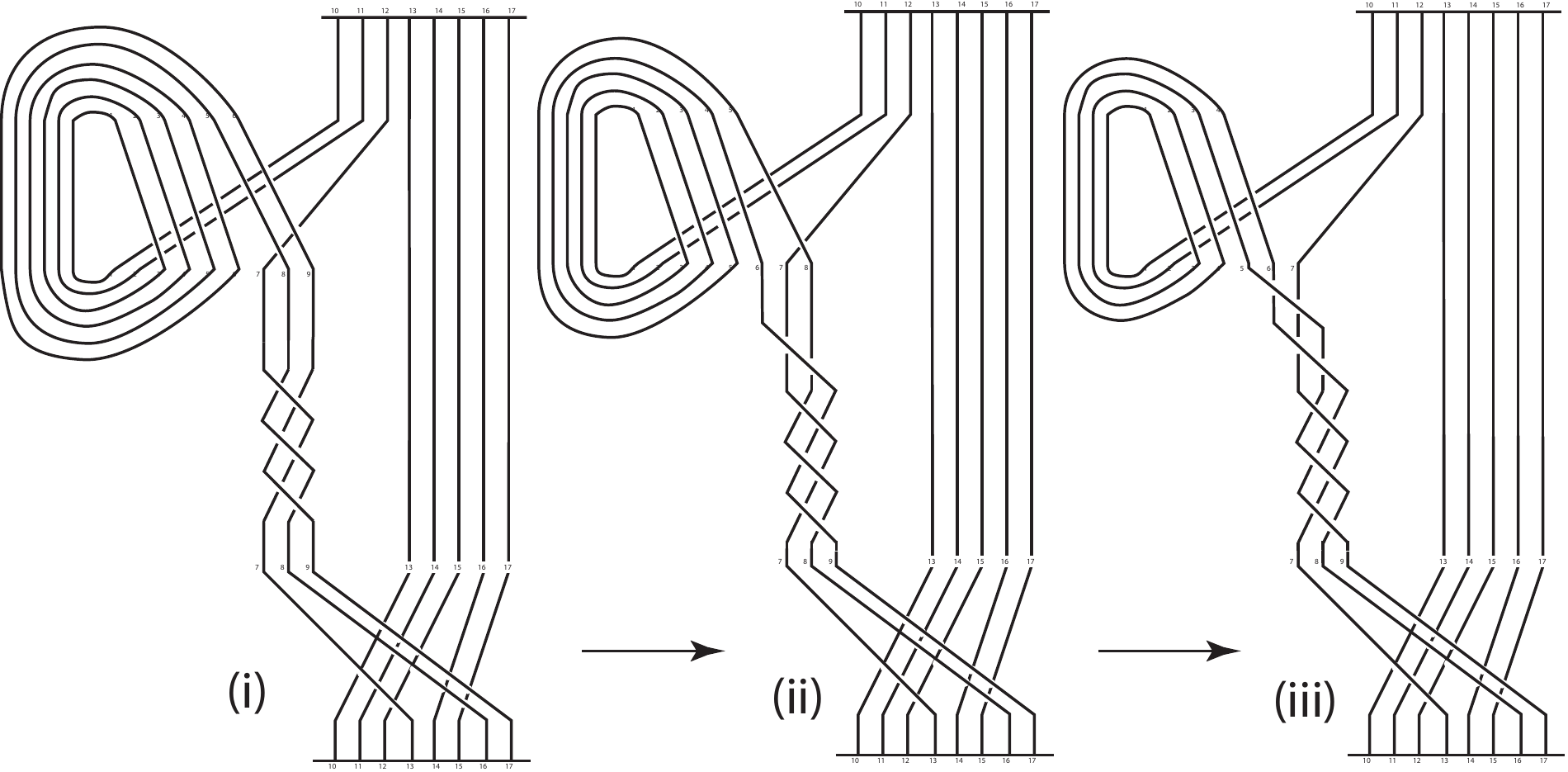}}
\caption{An example to show how a Lorenz link changes when we cut open the template}
\label{JBuncoil} 
\end{figure}

%%%%%%%%%%%%
Cut open the template along an orbit, as in the leftmost sketch in Figure~\ref{JBcutting-open}, splitting it into 4 pieces: the first carrying the LL strands,  the second the LR strands, the third  the  RL strands and the fourth the RR strands.   Stretch the cut-open template as in sketch (ii) in 
Figure~\ref{JBcutting-open}.   After the stretching we can see that the LR band can be separated from the LL band and `uncoiled', as in the passage from sketch (ii) to sketch (iii).   Sketch (iv)  shows what happens when the strands of the LL band are uncoiled, one strand at a time.  (This is illustrated by the example in Figure~\ref{JBuncoil}).  Finally, in sketch (v) of Figure~\ref{JBcutting-open} the RR and RL bands are cut open (as they were in Figure~\ref{JBbraidtemplate}).      A new braid structure emerges.  Observe that the LR braid is always one full twist of the LR band, so it depends only on the number of LR strands in a Lorenz knot.   

%%%%%%%%%%
\noindent Let's look at an example.  The braids that are supported on the template all have the same LR braid, but their LL-braids differ.
%%%%%%%%%%

\begin{figure}[htpb!]
\centerline{\includegraphics[scale=.60] {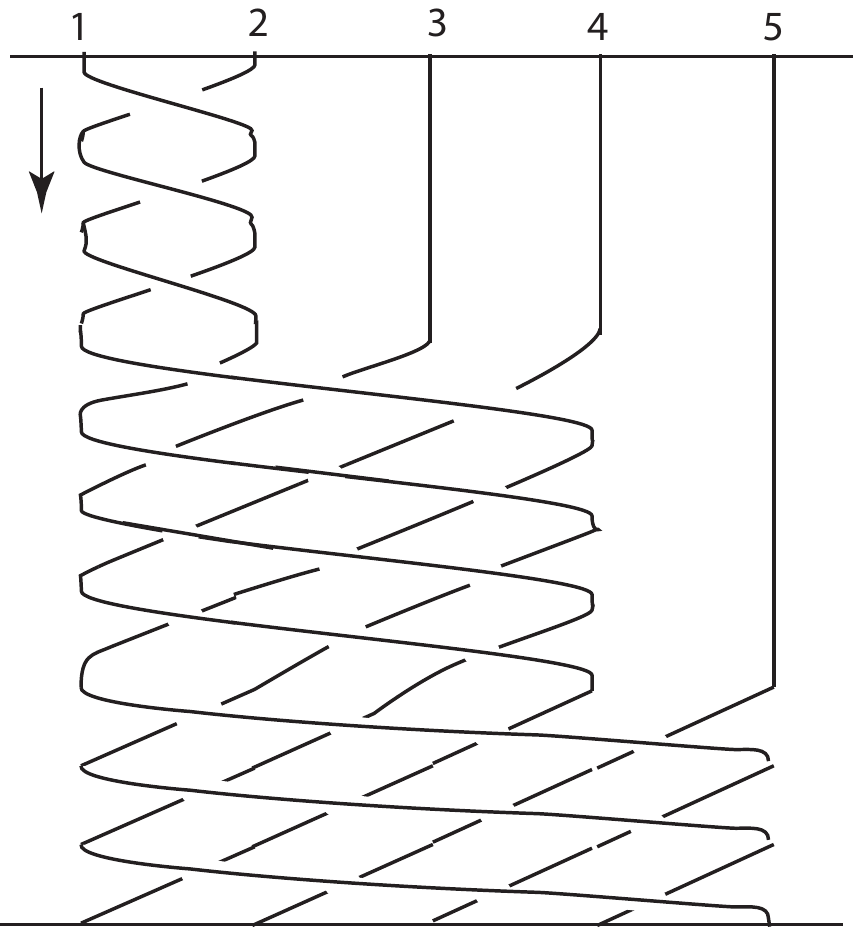}}
\caption{The T-braid with parameters ((2,3),(4,4),(5,3))} 
\label{JBT-braid}
\end{figure}

Ilya Kofman and the author have proved \cite{JBBKof09} that after the cutting and uncoiling processes that is illustrated in Figure~\ref{JBcutting-open} we obtain a {\it T-braid}.  Even more,   Lorenz links coincide, as a class, with links that we call {\it T-links}.  They are generalizations of torus links.  A torus link can be characterized, by a pair of integers $(p,q)$.  The most general T-link is defined by a sequence of integers $((p_1,q_1), (p_2,q_2),\dots, (p_r,q_r))$, where each $p_i < p_{i+1}$.   See Figure~\ref{JBT-braid} for an example when the parameters are ((2,3),(4,4),(5,3)).  This braid has 5 strands.
The T-link is the closure of a braid on $p_r$ strands.  First use the leftmost $p_1 < p_r$ strands ($p_1=2$ in the example) to form a torus braid of type $(p_1,q_1)$, then follow it by a torus braid of type $(p_2,q_2)$ on the leftmost $p_2$ strands, and so forth until, at the last step, all $p_r$ strands are used to form a torus braid of type $(p_r,q_r)$ .   The braid word (using the standard braid generators $\sigma_1,\dots,\sigma_{4}$ for $B_5$ that describes it is
$(\sigma_1)^3(\sigma_1\sigma_2\sigma_3)^4(\sigma_1\sigma_2\sigma_3\sigma_4)^3$ in the example.    In the general case, when the parameters are $(p_1,q_1),(p_2,q_2),\dots,(p_k,q_k))$, the braid will be:
$$(\sigma_1\sigma_2\dots,\sigma_{p_1-1})^{q_1}(\sigma_1\sigma_2\dots,\sigma_{p_2-1})^{q_2}\dots(\sigma_1\sigma_2\dots,\sigma_{p_k-1})^{q_k}$$
Braids with this general form are called {\it T-braids}, and their closures are called {\it T-links}. It is proved in \cite{JBBKof09}that  Lorenz links are in one-to-one-correspondence with T-links.  This is our third parametrization.

%%%%%%%%%%%%%%%%%%%%%%%%%%%%%%%%%%%%%%%%%%%%%%%%%%%%%%%%%%%%%%%%%%%%%%%

\section  {Applications}
\label{JBapplications}  
In this section we discuss some of the many consequences of our three viewpoints.

%%%%%%%%%%%%%
\subsection{The recognition problem for Lorenz knots}
In $\S$\ref{JBSS:questions} we asked how rich the class of Lorenz knots might be, and whether we could solve the recognition problem.  We give partial answers to both questions now. To explain what we know, we need several standard definitions:
%%%%%
\begin{itemize}
\item Let $c$ be the number of double points in a projection of a knot $K$.  The {\it crossing number} $c_{min}(K) $ is the smallest number of double points, over all possible projections of a knot $K$.  As we explained earlier  $c_{min}(K)$ is, by definition, a knot type invariant, but  it is in general very difficult to compute.
\item Every link can be represented as a closed braid, in many ways  \cite{JBRo}.  The {\it braid index} $n_{min}(K)$ is the smallest integer such that a knot or link $K$ can be represented as a closed braid on $n$
strands.   It too is by definition an invariant, and it too is in general difficult to compute.

\item Every knot is the boundary of a compact oriented surface. The {\it genus} $g(K)$ of a knot  $K$ is the genus of the surface $\Sigma$ of smallest genus such that $K =\partial \Sigma$, and again it's an invariant.
%%%%%%%%%%
\item A knot $K\subset \mathbb S^3$ is {\it fibered} if its complementary space $\mathbb S^3\setminus K$ has the structure of a surface bundle. If so, the fiber is a surface of minimum genus.  We note that there are knots that are not fibered, in fact the property of being fibered is quite special. 
%%%%%
\item  A knot is a {\it positive closed braid} if it has a closed braid projection in which all crossings are positive.   We note that every positive closed braid is fibered, and also that  there are fibered knots that are not positive closed braids. 
%%%%%
\end{itemize}
%%%%%%%%%%%
%%%%%%%%%%%
We turn to consequences for Lorenz knots. (To simplify things, we do not state the analogues for Lorenz links):
\begin{enumerate}
%%%
\item  The fact that Lorenz knots are precisely the knots that have natural embeddings in the Lorenz template shows that {\it Lorenz knots are positive closed braids}.  When we began to write this review it was not clear to us whether  every knot that is a positive closed braid was  a Lorenz knot,  so we asked some questions.  We learned from S. Jablon that the only 10-crossing knot that is Lorenz is the knot $10_{124}$.  Since the standard projection of the knot $10_{152}$ exhibits it as a positive closed braid, we conclude that:
\begin{itemize}
\item All knots$\supset$ Fibered knots $\supset$ Closed + braids $\supset$  Lorenz knots.
\end{itemize}
%%%
\item A consequence of positivity  is that for any closed positive $n$-braid representation of a Lorenz knot, the genus $g$, the crossing number $c$ of the given $n$-braid representative, and the braid index $n$ of the representative are related by the formula 
\begin{equation} \label{JBg-c-n formula}
2g(K)  =  n  +1 - c .
\end{equation}
Observe that $g(K)$ is a knot type invariant, but $n$ and $c$ are not.  From this it follows that  the genus $g(K)$ of a Lorenz knot $K$ can be computed from any projection in which all crossing are positive, e.g. the Lorenz template projection, by a simple counting argument.   In particular, using the first braid parametrization of a Lorenz knot (see Figure~\ref{JBLorenzbraid}), equation (\ref{JBg-c-n formula}) shows that:
\begin{equation}
2g(K) = |L|+|R| + 1 - \sum_{i=1}^{i=k} q_i(p_i-1).
\end{equation}
\item The integers $n_{min}(K)$ and $c_{min}(K)$ are (by definition) knot type invariants, but for arbitrary knots they are difficult to compute.  However, in the special case of Lorenz knots it is a consequence of the work done by Birman and Williams  in \cite{JBBiWil83} and by Franks and Williams in \cite{JBFrWil87}  that:
\begin{equation} \label{JBbraid-index}
n_{min}(K) = {\rm min} (|{\rm  LR}|,  |{\rm RL}|). 
\end{equation}
Thus we can compute  $n_{min}(K)$, too, from the Lorenz braid.   
\item In view of this, one might wonder how the invariant $c_{min}(K)$ is related to $g(K)$ and $n_{min}(K)$?
In general, the minimum crossing number of a knot is {\it not} realized at minimum braid index, however for Lorenz knots, it follows from Proposition 7.4 of \cite{JBMur91} that it is.  So we have $c_{min}(K) =  n_{min}(K) + 1 -2g(K)$.   
 Thus we can compute $c_{min}(K)$ too from the first parametrization of a Lorenz knot:
 \begin{equation}
 c_{min}(K) = \sum_{i=1}^k q_i(p_i-1) - |{\rm LL}| - |{\rm RR}|
 \end{equation}
 \end{enumerate}
\noindent Going beyond positivity, we know a little bit more:
\begin{enumerate}
\item [(5)] Lorenz knots have yet one further property that generalizes the fact that torus knots are types $(p,q)$ and $(q,p)$ have the same knot type, namely Lorenz knots have a very special symmetry of order 2, induced by rotating the template 180 degrees about the $z$ axis in  Figure~\ref{JBbraidtemplate}.  
\item [(6)] Pierre Dehornoy has investigated the zeros of the Alexander polynomial of a Lorenz knots.  See his arXiv preprint  \cite{JPDehor11} .
\end{enumerate}
%%%%%%
It seems possible that, with these two additional properties we already know  how to solve the recognition problem for Lorenz knots.   To give some evidence, we discuss a remarkable calculation that was done by Pierre Dehornoy and Etienne Ghys, with the help of Slavik Jablon.  They proved:
 \begin{itemize}
\item  Among the 1,701,936 prime knots having projections with $\leq 16$ crossings, only 19 are Lorenz, and among those only 7 are not torus knots!  In particular, they solved the recognition problem, when the crossing number is not  $\leq$16.   
\end{itemize}
However this may not say very much about the recognition problem because Lorenz knots generally have high crossing numbers.  For example, it was proved in \cite{JBBiWil83} that all algebraic knots are Lorenz, however there is not a single algebraic knot which is not a torus knot and has $c_{min} \leq 16$.

%%%%%%%%%%%%%%%
\subsection{The volume of $\mathbb S^3\setminus \mathcal K$ and the monodromy of $\mathcal K$} \label{JBss:hyperbolic volume}  
%%%%%%%%%%
Up to now, we have been discussing the topology of Lorenz knots.  This means that if $K$ is a Lorenz knot, we have been studying properties of the 3-manifold $S^3\setminus K$  that are independent of any particular choice of a metric on the manifold.  But in the early 1980's, William P. Thurston proved a remarkable theorem about the possibilities for assigning a unique, complete finite-volume metric to the complement of a hyperbolic  knot in $S^3$, also he proves that `most knots are hyperbolic'.  See \cite{JBScott} for a good introduction to this major topic, noting that what he calls the {\it Geometrization Conjecture} is now a theorem.  For hyperbolic knots the hyperbolic volume is a knot type invariant, which becomes an invariant for all knots by defining the volume to be zero when the knot is not hyperbolic.   A natural question, then, is to compute the volumes of hyperbolic Lorenz knots.  

The study of hyperbolic knot
complements has been a focal point for much recent work in 3-manifold
topology.  The very question changes the focus of
knot theory from the properties of diagrams to the geometry of the
complementary space.  Ideal tetrahedra are the natural building blocks
for constructing hyperbolic 3-manifolds, and ideal triangulations can
be studied with the help of computer programs.  In particular,  it has been learned that there are precisely 6,075
noncompact hyperbolic 3-manifolds that can be obtained by gluing the
faces of at most seven ideal tetrahedra \cite{JBCaHiWe99}, with recent extensions to  18,921 examples with eight tetrahedra by Culler and Dunflied.  For a hyperbolic knot, the
minimum number of ideal tetrahedra required to construct its
complement is a natural measure of its geometric complexity. It is shown in Table 1 of \cite{JBBKof09} that:
 
$\bullet$ Of the 201 hyperbolic knots in the \cite{JBCaHiWe99} census, at least 107 are Lorenz knots.  

\noindent The number 107 could be too small because, among the remaining 94
knots, it was not possible to decide whether five of them are or are not
Lorenz.  We remark that the known diagrams for the knots in question did not in any way
suggest the Lorenz template.  

There is an intuitive reason for this extraordinary data.  Torus knots, as observed earlier, have (by definition) hyperbolic volume 0.  One expects, then, that one way to construct hyperbolic knots of small volume is to `tweak' the diagram of a torus knot just a little bit.  Now recall that, using the first and third parametrizations,  Lorenz knots that are torus knots have parameters $(p,q)$, wherer $p$ and $q$ are coprime integers and $p< q$.  Therefore one might expect that Lorenz knots which have parameters $((1,2),(p_2,q_2))$ and are hyperbolic have very small volume, and indeed that is what the data in   \cite{JBBKof09} suggests.  This is, as we write,  a very active and interesting area of investigation. 

There is another aspect of the topology and geometry of Lorenz knots that is related but different from its volume, and we mention it very briefly.  As we observed, Lorenz knots are always fibered.   We then have the monodromy map  to study.  By a theorem of Thurston, a fibered knot $\mathcal K$ has a hyperbolic complement if and only if its monodromy map and all of its powers, acting on the fundamental group $\mathcal G = \pi_1(\Sigma)$ of the fiber, does not fix any conjugacy class in $\pi_1(\mathcal G)$.  The entropy of the monodromy map is then a new invariant of $\mathcal K$.  This topic is a big one, and space considerations prevent us from saying more, except to remark that just as Lorenz knots appear often in the census of fibered 3-manifolds having small but non-zero volume, so their monodromy maps appear often  in the study of surface diffeomorphisms having small positive entropy.   

%%%%%%%%%%%%%%%%%%%%%%%%%%%%%%%%%%%%%%%%%%%%%%%%%%%%%%%%%%%%%%%%%%%%%%%%%%%%%%%%%%%%%%%%%%%%%%%%%%%%%%%%%%%%%%%%%%%%%%%%%%%%%%%%%%%%%%%%%%%%%%%%%%%

\subsection{Modular knots and Lorenz knots}  

Modular knots were introduced by Etienne Ghys in \cite{JBGhy07} and \cite{JBGhy11}.  As we shall see, they are Lorenz knots in disguise.   In this section we explain his work, briefly.  

Recall that a model for hyperbolic 2-space $\mathbb H$ is $\{z = x+iy \in \mathbb C: y > 0 \}$ with the hyperbolic metric $ds = (\frac{1}{y})\sqrt {dx^2 + dy^2}.$  With this metric, a  {\it geodesic} in $\mathbb H$ is either a vertical ray in the upper half plane or a half-circle orthogonal to the $x$ axis.  Orientation-preserving isometries of $\mathbb H^2$ may be identified with $PSL(2,\mathbb Z)$, and the quotient space  $\mathbb H^2/PSL(2,\mathbb Z)$   is the {\it modular surface} $M$.    It has a natural metric coming from the metric on its covering space $\mathbb H^2$.

The  {\it geodesic flow on $M$} is defined by choosing a matrix $P  = \begin{bmatrix} a & b \\ c & d \end{bmatrix} \in PSL(2,\mathbb R)$ and sending $P \to H_tP,$   where  $H_t = \begin{bmatrix} e^t & 0 \\ 0 & e^{-t} \end{bmatrix}.$   \  
%Note that $P$ is hyperbolic precisely when $|a+d|>2$.  Projecting to $M$, there is then an induced  flow $\phi^t: M\to M$.  
This is the {\it geodesic flow $\phi^t$ on the modular surface $M$}.   An element $P$ defines a periodic orbit in the geodesic flow precisely when $H_t P = PA$ for some  hyperbolic matrix $A \in $PSL$(2,\mathbb Z)$.   Thus closed orbits in $\phi^t$ are defined by diagonal matrices $PAP^{-1}$, where $P\in$PSL$(2,\mathbb R)$, $A \in $PSL$(2,\mathbb Z)$.   The condition that $A$ be hyperbolic means that its trace in $>2$. 

In  \cite{JBGhy07} Etienne Ghys had the idea to lift closed orbits in the geodesic flow $\phi^t$ on $M$ to a related flow $\Phi^t$ on the unit tangent bundle $\tilde{M}$ of $M$,  a 3-manifold.  He used the known fact ($\S$ 3.1 of \cite{JBGhy07}) that $\tilde{M}$ is naturally isomorphic to $PSL(2,\mathbb R)/PSL(2,\mathbb Z)$, which is in turn known to be isomorphic to $\mathbb S^3 \setminus \mathcal T$, where $\mathcal T$ is the trefoil knot, i.e. the type (2,3) torus knot, embedded in $\mathbb S^3$ in a known canonical way.    This enabled him to study the closed orbits in $M$, a 2-manifold,  as knots in $\tilde{M} = \mathbb S^3 \setminus \mathcal T$.

\smallskip

\noindent {\bf Definition \cite{JBGhy07}:}  A {\it modular knot} is a closed orbit in the flow $\Phi^t$ on the unit tangent bundle $\tilde{M}\cong (\mathbb S^3\setminus \mathcal T)$ of the modular surface $M$.  

\smallskip

\noindent The geodesic flow on $M$ has been studied extensively by number theorists, and is of great interest.   However, while the topology and geometry that we just described is well-known, it seems as if nobody had ever really thought 
of its closed orbits as knots.  In particular, Ghys defines and studies the {\it Rademacher function}, an integer-valued function on closed orbits in the modular flow, via the geometry that we just described, demonstrating that its value on a closed orbit $K$ is in fact the linking number of $K$ with  the `missing trefoil' $\mathcal T$.

Lorenz knots have not entered into this picture yet, but we now show that their appearance is very natural and concrete. Modular knots are simple closed curves in $\tilde{M}  = \mathbb S^3 \setminus \mathcal T$.  Recall that the fundamental group of  $\mathbb S^3 \setminus \mathcal T$, admits the presentation 
$$ G = \pi_1(S^3 \setminus \mathcal T): < U,V; \ \ U^2 = V^3>.$$
The group $G$ is the free product of the cyclic groups generated by $U,V$ amalgamated along $C = U^2=V^3$,  where  $C$ generates the center of $G$.  Thus every free homotopy class in $G$ is represented by a cyclic word $W$ of the form
$C^k U V^{\epsilon_1}U V^{\epsilon_r} \cdots U V^{\epsilon_r}$, where $ \epsilon_i = \pm 1.$    

We have already shown that the choice of a modular knot corresponds to a choice of a hyperbolic matrix $A$ in the group $PSL(2,\mathbb Z)$, i.e. in the image of $G$ under the homomorphism $G \to PSL(2,\mathbb Z)$ whose kernel is the central element $C$.  That is, the choice of a cyclic word $\pm U V^{\epsilon_1}U V^{\epsilon_r} \cdots U V^{\epsilon_r}$, where $ \epsilon_i = \pm 1.$    Let $L = UV, \ \ R = UV^{-1}$.   Then  our cyclic word  goes over to a cyclic word in $L $ and $R$.   But then, by the parametrization given in $\S$\ref{JBsecond parametrization} above, each closed orbit  $K$ in the modular flow, lifted to $\tilde{M}$, determines a Lorenz knot!   

\smallskip

\noindent {\bf Theorem} (E. Ghys, \cite{JBGhy07}).   There is a one-to-one correspondence between modular knots and Lorenz  knots.    

\smallskip

\noindent There is more than this. Ghys proves that the Lorenz template occurs in the setting of modular knots, embedded in a natural way in $\mathbb S^3 \setminus \mathcal T$, and there is a natural way (using $L$-$R$ words) to embed the family of all modular knots, disjointly and simultaneously, onto the template.  

\begin{figure}[htpb!]
\centerline{\includegraphics[scale=.70]{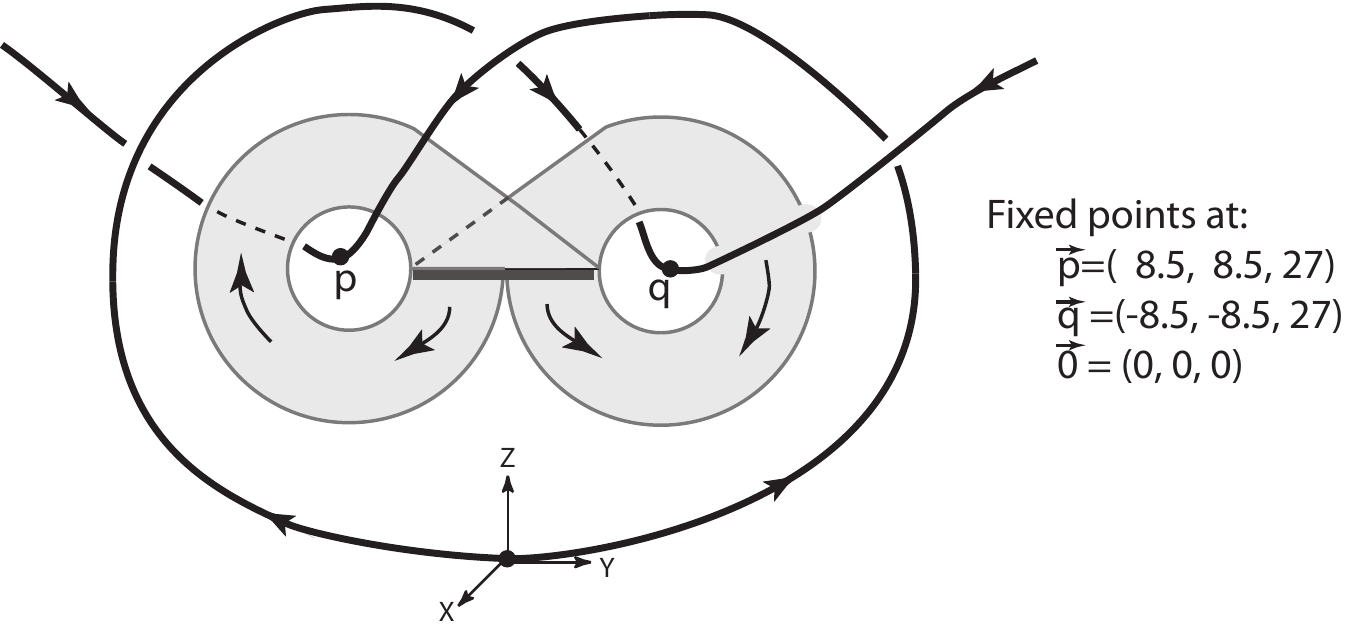}}
\caption{The missing treoil.} \label{JBmissing-trefoil}
\end{figure}

There are two issues about this picture that remain somewhat mysterious.  The first is that (a) Ghys does not address the issue that, if they coincide, then there ought to be a `missing trefoil' in the Lorenz flow.  The second is that  (b) Ghys does not actually prove that the Lorenz flow in $\mathbb S^3$ and the modular flow coincide.  Regarding (a), Tali Pinsky \cite{JBPi10} has computer evidence that in fact there is a missing trefoil very naturally embedded in $\mathbb S^3$ with respect to the Lorenz flow. In \cite{JBPi10} she conjectures that the invariant curves connecting the three fixed points of the Lorenz flow (one in the center of each `ear' of the template, and one on the axis of symmetry) join up to form such a trefoil.  See Figure~\ref{JBmissing-trefoil}. Since the curves in question are invariant curves, closed orbits in the flow would necessarily avoid them.  We look forward to her proof that the union of the invariant curves is actually a knot.  
As for (b), while Ghys is a bit vague about how he produced the beautiful pictures in \cite{JBGhy11}, the picture proofs `say more than a thousand words'.  We recommend them strongly to any reader with a small interest in this topic, as a very non-traditional introduction to Lorenz knots.  

In a very different direction, his discoveries reveal a new fact about the Rademacher function: As mentioned very early in this review, the Lorenz flow, which we now reinterpret as the modular flow on $\tilde{M} \cong \mathbb S^3 \setminus \mathcal T$ is chaotic in the sense that it has an extremely sensitive dependence on initial conditions.  This is true for its closed orbits too, so that if two points are very close on the modular surface, and both lie on closed orbits, their linking numbers with $\mathcal T$ might be very different.  Thus it's a canonical example of what one might call a {\it chaotic function}.

\subsection{Generalizations of Lorenz knots}  In $\S$\ref{JBSchaos} of this paper we introduced templates, via the example of Lorenz knots.  We never even hinted that there were templates for  flows on $\mathbb S^3$ different from the Lorenz flow.  But in fact the Lorenz template is not an isolated tool.  We give several examples.

(i) The paper \cite{JBBiWil83}, which initiated the study of Lorenz knots, had a twin \cite{JBBiWil-II-83} in which a different flow  on $S^3$ and its associated template and family of closed orbits was discussed and studied.   The second flow arose from a phenomenon which we have already encountered.  The knot $\mathcal T' = 4_1$ (see Figure~\ref{JBknottables}) is a fibered knot of genus 1, so the fiber is a once-punctured torus.   The monodromy map, lifted to the universal cover of the closed torus, i.e. the Eucldean plane $\mathbb R^2$, is the linear map $ \begin{bmatrix} 2 & 1 \\ 1 & 1 \end{bmatrix}$.   This map has a dense subset of periodic points,  from which it follows that the closed orbits in the associated flow on $(\mathbb S^3 \setminus \mathcal T')$, defined by pushing the fiber around the knot, are dense in the flow. Thus, just as Ghys had shown that the modular flow lifted to $(\mathbb S^3\setminus \mathcal T)$,  the flow associated to the fibration is a flow on $(\mathbb S^3 \setminus \mathcal T')$.   In \cite{JBBiWil83} the author and Williams constructed a template for the later.  Subsequently,  R. Ghrist proved that it includes {\it all} knots and links.  Its associated template is said to be {\it universal}, and it seems likely that generic templates are, in fact, universal.  

(ii) The research monograph \cite{JBGhHS97} is dedicated to the general topic of knots and links that are defined by templates, through 1997.   Referenced in that monograph are papers that deal with the problem of passing from a particular template  to an associated flow, a non-trivial question.   There are a host of open problems, as we write, about templates that are universal and those that are not.

(iii) In a related but different direction, Pinsky has investigated, in \cite{JBPi10}, templates for geodesic flows on certain hyperbolic manifolds which are different from the geodesic flow on the modular surface, and she has found templates there too, and also analogues of the `missing trefoil' in the unit tangent bundle of the modular surface.  There is a big world out there, and a great deal of structure, waiting to be discovered! 

\section{Acknowledgements}  Thanks go to to Ilya Kofman for reading a very preliminary draft of this manuscript  and commenting on it, and for stimulating discussions; to Morwen Thistlethwaite for the knots in Figure~\ref{JBsameJonespoly}; to an unknown colleague for the picture in Figure~\ref{JBbeautifulpic}, which we downloaded from the internet.

%%%%%%%%%%%%%%%%%%

\end{document}